\documentclass[12pt,letterpaper,reqno]{amsart}
\usepackage{amscd}
\usepackage{latexsym}
\usepackage{amsfonts}
\usepackage{amssymb}
\newtheorem{theorem}{Theorem}
\newtheorem{proposition}[theorem]{Proposition}
\newtheorem{lemma}[theorem]{Lemma}
\newtheorem{corollary}[theorem]{Corollary}
\theoremstyle{definition}
\newtheorem{definition}[theorem]{Definition}

\newtheorem{remark}[theorem]{Remark}
\newcommand{\rk}{\mathrm{rk}}

\newcommand{\Hom}{\mathrm{Hom}}
\newcommand{\Inn}{\mathrm{Inn}}
\newcommand{\ind}{\mathrm{Ind}}

\newcommand{\tr}{\mathrm{tr}}

\newcommand{\aut}{\mathrm{Aut}}

\newcommand{\Rep}{\mathrm{Rep}}

\newcommand{\Char}{\mathrm{Char}}
\newcommand{\Syl}{\mathrm{Syl}}
\begin{document}
\title[$\textrm{Qd}(p)$-free rank two groups have homotopy rank two]
{$\textrm{Qd}(p)$-free rank two finite groups act freely on a homotopy product of two spheres}
\author{Michael A. Jackson}
\address{Department of Mathematics, Hylan Building, University of Rochester, Rochester, NY 14627, USA\\
mjackson@math.rochester.edu}

\begin{abstract}
A classic result of Swan states that a finite group $G$ acts freely on a finite homotopy sphere if and 
only if every abelian subgroup of $G$ is cyclic. Following this result, Benson and Carlson 
conjectured that a finite group $G$ acts freely on a finite complex with the homotopy type of $n$ spheres if 
the rank of $G$ is less than or equal to $n$. Recently, Adem and Smith have shown that 
every rank two finite $p$-group acts freely on a finite complex with the homotopy type of two spheres. In this paper 
we will make further progress, showing that rank two groups that are $\textrm{Qd}(p)$-free act freely on a finite homotopy
product of two spheres. 

{\it 2000} MSC: 57Q91, 55R25, 20C15, 20E15
\end{abstract}
\maketitle

\section{Introduction}
Recall that the $p$-rank of a finite group $G$, $\rk _{p}(G)$, is the largest rank of an elementary 
abelian $p$-subgroup of $G$ and that the rank of a finite group $G$, $\rk (G)$, is the maximum of 
$\rk _{p}(G)$ taken over all primes $p$. 
We can define the homotopy rank of a finite group $G$, $\textrm{h} (G)$, to be the minimal 
integer $k$ such that $G$ acts freely on a finite CW-complex 
$Y \simeq \mathbb{S}^{n_1} \times \mathbb{S}^{n_2}\times \cdots \times \mathbb{S}^{n_k}$. 
Benson and Carlson \cite{bc:cmc} state a conjecture that for any finite group 
$G$, $\rk (G)=\textrm{h}(G)$. The case of this conjecture when $G$ is a rank one group
is a direct result of Swan's theorem \cite{s:prfg}. Benson and Carlson's conjecture has also been 
verified by Adem and Smith \cite{as:ospc,as:pcga} in the case when $G$ is a rank two $p$-group. 
In addition, recall that Heller \cite{h:nso} has shown that for a finite group $G$, if $\textrm{h}(G)=2$, 
then $\rk (G)=2$. In this paper we will verify this same conjecture for all rank two groups that 
do not contain a particular type of subquotient.\\
\\
Before we state our main theorem, we need two definitions.

\begin{definition}
A \emph{subquotient} of a group $G$ is a factor group $H/K$ where $H,K \subseteq G$ with $K\lhd H$. A subgroup 
$L$ is said to be involved in $G$ if $L$ is isomorphic to a subquotient of $G$. In particular, for a prime $p$, we say that 
$L$ is $p'$-involved in $G$ if $L$ is isomorphic to a subquotient $H/K$ of $G$ where $K$ has order relatively prime to $p$. 
\end{definition}

\begin{definition}\label{def:qdp}
Let $p$ be a prime and $E$ be a rank two elementary abelian $p$-group. $E$ is then also a two-dimensional vector space over 
$\mathbb{F}_p$. We define the quadratic group $\textrm{Qd}(p)$ to be the semidirect product of $E$ by the special linear group 
$SL_2(\mathbb{F}_p)$.
\end{definition}

A classic result of Glauberman shows that if for some odd prime $p$ $\textrm{Qd}(p)$ is not involved in a finite group 
$G$, then the $p$-fusion is controlled by the normalizer of a characteristic $p$-subgroup \cite{g:psg,g:scps}. Although 
we do not directly use Glauberman's result here, there are some connections as $p$-fusion plays a integral role in our 
investigation (see Remark \ref{rem:glaub}). We can now state our main theorem.
\begin{theorem}\label{thm:main}
Let $G$ be a finite group such that $\rk (G)=2$. $G$ acts freely on a finite CW-complex 
$Y \simeq \mathbb{S}^{n_1} \times \mathbb{S}^{n_2}$ unless for some odd prime p, 
$\textrm{Qd}(p)$ is $p'$-involved in $G$.
\end{theorem}

The following corollary, which can also be shown using more direct methods, follows immediately from 
Theorem \ref{thm:main}:

\begin{corollary}
If $G$ is a finite group of odd order with $\rk (G)=2$, then $G$ acts freely on a finite CW-complex 
$Y \simeq \mathbb{S}^{n_1} \times \mathbb{S}^{n_2}$.
\end{corollary}

I am grateful to Radha Kessar for her help in clarifying the argument in Section \ref{sec:odd} as well as pointing 
out some relevant connections to group theory. I would also like to thank Antonio Viruel and Antonio Diaz for helpful comments, 
which aided in the completion of the proofs in Section \ref{sec:odd}.

\section{Preliminaries}\label{sec:pre}
Recall that Adem and Smith have shown that if $G$ is a rank two $p$-group, then $h(G)=2$ \cite{as:ospc,as:pcga}. 
In their proof they use the following theorem: 

\begin{theorem}[Adem and Smith \cite{as:pcga}]\label{thm:as1}
Let $G$ be a finite group and let $X$ be a finitely dominated, simply connected $G$-CW complex such that 
every isotropy subgroup has rank one. Then for some integer $N>0$ there exists a finite CW-complex 
$Y\simeq \mathbb{S}^N \times X$ and a free action of $G$ on $Y$ such that the projection $Y\rightarrow X$ 
is $G$-equivariant.
\end{theorem}

In the present work, we apply Theorem \ref{thm:as1} for a rank two group by finding an 
action of the group on a 
CW-complex $X\simeq S^n$ such that the isotropy subgroups are of rank one. Adem and 
Smith in \cite{as:pcga} also developed a sufficient condition for this type of action to exist.\\
\\
To explore this sufficient condition, we start by letting $G$ be a finite group. Given a 
$G$-CW complex $X \simeq \mathbb{S}^{N-1}$, there is 
a fibration
$X \rightarrow X\times _{G} EG \rightarrow BG$. This fibration gives an Euler class 
$\beta $, which is a cohomology class in $H^{N}(G)$. This cohomology class is the 
transgression in the Serre spectral sequence of the fundamental class of the fiber. 
Regarding this class, Adem and Smith made the following definition:
\begin{definition}[Adem and Smith \cite{as:pcga}]
Let $G$ be a finite group. The Euler class $\beta$ of a $G$-CW complex
$X \simeq \mathbb{S}^{N-1}$ is called effective if the Krull dimension
of $H^{\ast}(X\times _{G} EG, \mathbb{F}_{p})$ is less than $\rk(G)$ for all $p\bigr| |G| \bigl.$.
\end{definition}

Adem and Smith have proven the following lemma characterizing effective Euler classes: 
\begin{lemma}[Adem and Smith \cite{as:pcga}]\label{lem:as2}
Let $\beta \in H^N (G)$ be the Euler class of an action on a finite dimensional $X \simeq 
\mathbb{S}^{N-1}$. $\beta $ is effective if and only if every maximal rank elementary abelian 
subgroup of $G$ acts without stationary points.
\end{lemma}

\begin{remark}\label{rmk:as2}
By applying Lemma \ref{lem:as2} and Theorem \ref{thm:as1}, we see that any finite group of rank two, 
having an effective Euler class of an action on a finite dimensional $X \simeq \mathbb{S}^{N-1}$, acts 
freely on a finite CW-complex $Y\simeq \mathbb{S}^{m}\times \mathbb{S}^{N-1}$.
\end{remark}

We now make a topological definition that will allow us to classify certain actions on homotopy spheres.

\begin{definition}
Consider the set of maps from a space $X$ to a space $Y$. Two such maps $f$ and $g$ are called homotopic 
if there exists a continuous function $H: X \times [0,1] \rightarrow Y$ such that for all $x\in X$, $H(x,0)=f(x)$ and 
$H(x,1)=g(x)$. This induces an equivalence relation on the set of maps from $X$ to $Y$, which makes two maps 
equivalent if they are homotopic. The resulting equivalence classes are called homotopy classes of maps from 
$X$ to $Y$; the set of such classes is denoted $[X,Y]$. 
\end{definition}

Recall that for a group $G$, $BG$ is the classifying space of $G$. 
Also $[BG,BU(n)]$ represents 
the homotopy classes of maps from $BG$ to $BU(n)$, which may also be written as 
$\pi _0 \Hom (BG,BU(n))$ (see \cite{dz:mbcs}). In addition, 
recall that $[BG, BU(n)]$ is isomorphic to the set of isomorphism classes of rank 
$n$ complex vector bundles over $BG$. In this discussion, complex vector bundles over 
$BG$ will be thought of as elements of $[BG,BU(n)]$.\\
\\
Recall that the Euler class of the universal bundle over $BU(n)$ is 
the top Chern class, which lies in $H^{2n}(BU(n))$. We will call this class $\xi $.
A map $\varphi : BG \rightarrow BU(n)$ induces a map $\varphi ^{\ast}: H^{2n}(BU(n))\rightarrow H^{2n}(BG)=H^{2n}(G)$. 
Thus the image of $\xi $ under $\varphi ^{\ast}$ is the Euler class of the vector bundle over 
$BG$ that corresponds to $\varphi $. Let $\beta _{\varphi} = \varphi ^{\ast}(\xi)$ be this 
class, which is the Euler class of a $G$-CW complex that is the unit sphere bundle of the 
vector bundle. This $G$-CW complex is the fiber as discussed at the outset of Section \ref{sec:pre}. 
Recall that if $\varphi $ and $\mu $ are homotopic maps from $BG$ to $BU(n)$, then 
$\beta _{\varphi}=\beta _{\mu }$; therefore, the Euler class of a homotopy class of maps is 
well defined.\\
\\
Recall that $G$ is a finite group, 
$p$ is a prime dividing the order of $G$, and $G_{p}$ is a Sylow $p$-subgroup of $G$. 
In addition we will use $\Char _{n}(G_{p})$ to represent the set of 
degree $n$ complex characters of $G_{p}$ and will let $\Char _{n}^{G}(G_{p})$ represent the 
subset consisting of those degree $n$ complex characters that are the 
restrictions of class functions of $G$. Also we use the notation $\Rep (G,U(n))$ for the degree $n$ 
unitary representations of the group $G$, which is given by the equality 
$\Rep (G,U(n))=\Hom (G, U(n))/\Inn (U(n))$. 
Recall that Dwyer and Zabrodsky \cite{dz:mbcs} have shown that for any $p$-subgroup $P$,
\[ \Char _n(P) \cong \Rep (P,U(n))\stackrel{\cong}{\rightarrow} [BP,BU(n)]\] 
where the map is $\rho \mapsto B\rho$.
So if $\varphi$ is a map $BG \rightarrow BU(n)$ and if $P$ is a $p$-subgroup of $G$ for some 
prime $p$, $\varphi |_{BP}$ is induced by a unitary representation of $P$. 
The following proposition gives the characterization of effectiveness for the Euler class 
$\beta _{\varphi}$.
\begin{proposition}[{\cite[section 4]{as:ospc}}]\label{prop:effec}
An Euler class $\beta_{\varphi} $ is effective if and only if for each elementary abelian 
subgroup $E$ of $G$ with $\rk(E)=\rk(G)$, there is a unitary representation $\lambda :E 
\rightarrow U(n)$ such that both $\varphi _E = B \lambda$ and $\lambda $ does not have 
the trivial representation as one of its irreducible constituents.
\end{proposition}

As we proceed we will want to work separately with each prime. 
Definition \ref{def:peff} will be useful in this endeavor.

\begin{definition}\label{def:peff}
Let $G$ be a finite group. The Euler class $\beta _{\varphi}$ is called p-effective if for each 
elementary abelian $p$-subgroup $E$ of $G$ with $\rk(E)=\rk(G)$ there is a unitary representation 
$\lambda :E \rightarrow U(n)$ such that both $\varphi _E = B \lambda$ and $\lambda $ does not have 
the trivial representation as one of its irreducible constituents.
\end{definition}

\begin{remark}
It is clear then that $\beta _{\varphi}$ is effective if and only if it is $p$-effective for 
each $p\bigr| |G| \bigl.$. Also if $\rk _{p}(G)<\rk(G)$, then any $\beta _{\varphi}$ is $p$-effective 
since it trivially satisfies Definition \ref{def:peff}; therefore, we can restate  
Proposition \ref{prop:effec}: $\beta _{\varphi}$ is effective if and only if it is $p$-effective for each 
prime $p$ such that $\rk _{p}(G)=\rk(G)$.
\end{remark}

Definition \ref{def:fus} concerns the characters of subgroups of $G$. 
\begin{definition}\label{def:fus}
Let $H\subseteq G$ be a subgroup. Let $\chi $ be  a character of $H$. We say that $\chi $
\emph{respects fusion in $G$} if for each pair of elements $h$, $k \in H$  that are conjugate in the group $G$, $\chi (h)=\chi (k)$.
\end{definition}

\begin{remark}\label{rem:fus}
Let $H\subseteq G$. A sufficient condition for a character of $H$ to respect fusion in 
$G$ is for the character to be constant on elements of the same order.
\end{remark}

Notice that for $\chi \in \Char _{n}(G_{p})$, $\chi $ respects fusion in $G$ if and only if 
$\chi \in \Char _{n}^{G}(G_{p})$. In a previous paper \cite{j:qsfg}, the author has 
shown the following theorem.

\begin{theorem}[Jackson {\cite[Theorem 1.3]{j:qsfg}}]\label{thm:surj}
If $G$ is a finite group that does not contain a rank three elementary abelian subgroup, 
then the natural mapping 
$$\psi _{G}:[BG,BU(n)] \rightarrow \prod _{p\bigr| |G| \bigl.}\Char _{n}^{G}(G_{p})$$ is a surjection.
\end{theorem}
Combining Theorem \ref{thm:surj} and Definition \ref{def:peff}, 
we get the following theorem.

\begin{theorem}\label{thm:peffimp}
Let $G$ be a finite group with $\rk (G)=2$, $p$ a prime dividing $|G|$, and $G_p$ a Sylow $p$-subgroup of $G$. 
If there is a character $\chi $ of $G_p$ that respects fusion in $G$ and has the property that 
$[\chi |_{E},1_{E}]=0$ for each $E\subseteq G_p$ elementary abelian with $\rk (E)=2$,
then there is a complex vector bundle over $BG$ whose Euler class is $p$-effective.
\end{theorem}
\begin{remark}\label{rmk:reg}
Let $p$ be a prime, $G_p$ a $p$-group of rank $n$, and $\chi $ a character of $G_p$. A sufficient condition 
to guarantee that $[\chi |_{E},1_{E}]=0$ for each $E\subseteq G_p$ elementary abelian with $\rk (E)=\rk(G_p)$ 
is that $\chi (1)=c (p^n-1)$ and $\chi (g)=-c $ for each element of $G_p$ of order $p$. In this case, $\chi |_E =
c(\rho -1)$ where $\rho $ is the regular character of the group $E$.
\end{remark}

Theorem \ref{thm:peffimp} leads to definition \ref{def:peffch}:
\begin{definition}\label{def:peffch}
Let $G$ be a finite group, $p$ a prime dividing $|G|$, and $G_p$ a Sylow $p$-subgroup of $G$. 
A character $\chi $ of $G_p$ is called a \emph{$p$-effective character of $G$} 
under two conditions: $\chi $ respects 
fusion in $G$; and for each $E\subseteq G_p$ elementary abelian with $\rk (E)=\rk(G)$, 
$[\chi |_{E},1_{E}]=0$.
\end{definition}

\begin{remark}\label{rem:peff}
Let $G$ be a finite group, $p$ a prime dividing $|G|$, and $G_p$ a Sylow $p$-subgroup of $G$. 
A character $\chi $ of $G_p$ is $p$-effective if $\chi $ is both constant on elements of the same order 
and $[\chi |_{E},1_{E}]=0$ for each $E\subseteq G_p$ elementary abelian with $\rk (E)=\rk(G)$ . 
This follows from Remark \ref{rem:fus}.
\end{remark}

Applying Definition \ref{def:peffch}, Theorem \ref{thm:peffimp} can be restated as follows: If $G$ has a 
$p$-effective character, then there is a complex vector bundle over $BG$ whose Euler class is 
$p$-effective.\\
\\
We have now essentially proved Theorem \ref{thm:peffch}, which is a 
modified form of a theorem by Adem and Smith \cite[Theorem 7.2]{as:pcga}.

\begin{theorem}[{See \cite[Theorem 7.2]{as:pcga}}]\label{thm:peffch}
Let $G$ be a finite group with $\rk (G)=2$. If for each prime $p$ dividing the order of $G$ there exists 
a $p$-effective character of $G$, then $G$ acts freely on a finite CW-complex 
$Y\simeq \mathbb{S}^{N_1}\times \mathbb{S}^{N_2}$.
\end{theorem}

\begin{proof} For each prime $p$ dividing the order of $G$, let $\chi _p$ be a $p$-effective character of $G$. 
There exists an integer $n$ such that $\chi _p (1)$ divides $n$ for each $p$ dividing the order of $G$ giving 
$$(\frac{n}{\chi _{p_1}(1)}\cdot \chi _{p_1}, \frac{n}{\chi _{p_2}(1)}\cdot \chi _{p_2}, \dots ,\frac{n}{\chi _{p_k}(1)}\cdot \chi _{p_k}
) \in {\prod _{p\bigr| |G| \bigl.}\Char _{n}^{G}(G_{p})}.$$
By Theorem \ref{thm:surj}, there exists an element $\varphi \in [BG, BU(n)]$ such that 
$$\psi _G(\varphi ) = (\frac{n}{\chi _{p_1}(1)}\cdot \chi _{p_1}, \frac{n}{\chi _{p_2}(1)}\cdot \chi _{p_2}, \dots ,\frac{n}{\chi _{p_k}(1)}\cdot \chi _{p_k}).$$
The Euler class of the homotopy class $\varphi $ is then an effective Euler class of $G$.
Theorem \ref{thm:peffch} then follows from Remark \ref{rmk:as2}. 
\end{proof}

Theorem \ref{thm:peffch} reduces the problem of showing that a rank two finite group $G$ acts 
freely on a finite CW-complex $Y\simeq \mathbb{S}^{N_1}\times \mathbb{S}^{N_2}$ to the problem 
of demonstrating that for each prime $p$ there is a $p$-effective character of $G$. 
It is important to state that not all rank two finite groups have 
a $p$-effective character for each prime $p$. An example of a finite group that does not 
contain a $p$-effective character for a particular prime is given in the following remark, which 
will be proven in Section \ref{sect:count}.
\begin{remark}\label{rem:psl3}
Let $p$ be an odd prime. The group $PSL_{3}(\mathbb{F}_{p})$ does not have a $p$-effective 
character.
\end{remark}

The groups $PSL_{3}(\mathbb{F}_{p})$, for odd primes $p$, are the only finite simple groups 
of rank 2 that were not shown by Adem and Smith to act freely on a finite CW-complex 
$Y\simeq \mathbb{S}^{N_1}\times \mathbb{S}^{N_2}$ \cite{as:ospc,as:pcga}.
\section{A sufficient condition for $p$-effective characters}\label{sec:suf}

In Section \ref{sec:suf}, we establish a sufficient condition for the 
existence of a $p$-effective character by first recalling the following definitions.

\begin{definition}\label{def:sc}
Let $G$ be a finite group and $H$ and $K$ subgroups such that $H\subset K$. We say that 
$H$ is strongly closed in $K$ with respect to $G$ if for each $g\in G$, 
$gHg^{-1} \cap K \subseteq H$.
\end{definition}

\begin{remark}\label{rmk:wc}
Let $G$ be a finite group and $H$ and $K$ subgroups such that $H\subset K$. We say that 
$H$ is weakly closed in $K$ with respect to $G$ if for each $g\in G$ with $gHg^{-1} \subseteq K$,  
$gHg^{-1}=H$. For a subgroup $H$ of prime order, $H$ is strongly closed in $K$ with respect to $G$ if and only if $H$ 
is weakly closed in $K$ with respect to $G$.
\end{remark}

\begin{definition}
Let $P$ be a $p$-group. We use $\Omega _1 (P)$ to denote a $p$-subgroup of $P$ generated by all of the order $p$ elements 
of $P$. In particular, if $P$ is abelian, $\Omega _1 (P)$ is elementary abelian, and if $P$ is cyclic, then $\Omega _1 (P)$ has 
order $p$.
\end{definition}

\begin{lemma}\label{lem:omega}
If $E$ is an elementary abelian subgroup of a $p$-group $P$, then  $\langle E,  \Omega _1 (Z(P))\rangle $
is an elementary abelian subgroup of $P$. In particular, if $E$ is a maximal elementary abelian 
subgroup of $P$ under inclusion, then $\Omega _1(Z(P)) \subseteq E$.
\end{lemma}

\begin{proof} Let $Z=\Omega _1 (Z(P))$. Since $Z$ is central, the subgroup $\langle E , \Omega _1 (Z(P))\rangle $ is abelian. 
The product is generated by elements of order $p$ and so must be elementary abelian. 
\end{proof}

\begin{corollary}\label{cor:omega}
Let $P$ be a $p$-group with $\rk(P)=\rk(Z(P))$. $\Omega _1(Z(P))$ is the unique elementary abelian 
subgroup of $P$ that is maximal under inclusion. In addition, $\Omega _1(P)=\Omega _1(Z(P))$, 
and $\Omega _1(Z(P))$ is strongly closed in $P$ with respect to $G$.
\end{corollary}

\begin{proof} If $E$ is a maximal elementary abelian subgroup of $P$ under inclusion, then $\Omega _1(Z(P)) 
\subseteq E$ by Lemma \ref{lem:omega}. Equality holds since $\rk (E) \leq \rk (P) = \rk(\Omega _1(Z(P)))$.
\end{proof}

Again we let $p$ be a prime dividing the order of a finite group $G$ and let $G_p\in \Syl _p (G)$. 
We show that if there is a 
subgroup of $Z(G_p)$ that is strongly closed in $G_p$ with respect to $G$, then $G$ has a $p$-effective
character in the following proposition.

\begin{proposition}\label{prop:scs}
Let $G$ be a finite group, $n=\rk(G)$, $p$ a prime divisor of $|G|$ with $\rk _{p}(G)=n$, 
and $G_p$ a Sylow $p$-subgroup of $G$. If there exists $H \subseteq Z(G_p)$ 
such that $H$ is non-trivial and strongly closed in $G_p$ with respect to $G$, 
then $G$ has a $p$-effective character.
\end{proposition}

\begin{proof}[Proof of Proposition \ref{prop:scs}] Let $H \subseteq Z(G_p)$ such that $H$ is non-empty and 
strongly closed in $G_p$ with respect to $G$. Let $K=\Omega _1 (H)$. 
So there is an $i\leq n$ with $K\cong (\mathbb{Z}/p\mathbb{Z})^{i}$. 
Let $\chi =\sum_{\psi \in Irr(K)\setminus \{1_{K}\}} \psi $; therefore, $[\chi , 1_{K}]=0$. Let 
$\varphi = \ind ^{G_p}_{K} \chi $.
Notice that for each $y\in K\setminus \{e\}$, $\chi (y)=-1$ \cite[4.2.7.ii]{g:fg}. 
In addition since $K\subseteq Z(G_p)$, $\varphi (y)=-[G_p:K]$  for each $y\in K\setminus \{e\}$
\cite[4.4.3.i]{g:fg}. $K$ is strongly closed in $G_p$ with respect to $G$ because $H$ is strongly closed 
in $G_p$ with respect to $G$. $\varphi $ respects fusion in $G$; 
therefore, for each $z\in G_p\setminus K$, $\varphi (z)=0$.\\
\\
Let $E$ be a subgroup of $G_p$ such that $E\cong (\mathbb{Z}/p\mathbb{Z})^{n}$.
By Lemma \ref{lem:omega}, $K\subseteq E$ since $K\subseteq \Omega _1(Z(G_p))$.
Now suppose $\exists \; E \subseteq G_p$ such that both $E\cong (\mathbb{Z}/p\mathbb{Z})^{n}$ and 
$[\varphi | _{E} , 1_{E}] > 0$. $K \subseteq \Omega _1(Z(G_p))$ implies that $K\subseteq E$; therefore, 
$[(\varphi | _{E}) |_{K}, 1_{K}] >0$. Note that 
$(\varphi | _{E}) |_{K}=\varphi | _{K}$ and that $\varphi | _{K} = [G_p:K] \chi$; thus 
$[\varphi |_{K}, 1_{K}] = [G_p:K] [\chi , 1_{K}] =0$.  This is a 
contradiction showing that if $E$ is an elementary abelian $p$-subgroup of $G_p$ of rank $n$, 
$[\varphi |_{E},1_{E}]=0$; therefore, $\varphi $ is a $p$-effective character of $G$. 
\end{proof}

We also prove a Corollary to Proposition \ref{prop:scs}. 

\begin{corollary}\label{cor:r2c}
Let $G$ be a finite group, $p$ a prime divisor of $|G|$ with $\rk (G)=\rk _{p}(G)=n$, 
and $G_p$ a Sylow $p$-subgroup of $G$. If $\rk (Z(G_p)) =n$, 
then $G$ has a $p$-effective character.
\end{corollary}

\begin{proof} By Corollary \ref{cor:omega}, $\Omega _1 (Z(G_p))$ is strongly closed in $G_p$ with respect to 
$G$. Corollary \ref{cor:r2c} then follows from Proposition \ref{prop:scs}.
\end{proof}

In Sections \ref{sec:odd} and \ref{sec:p2}, we look at the possible types of rank two groups where 
this sufficient condition holds. In Section \ref{sec:odd} we look at the case where $p$ 
is an odd prime, while in Section \ref{sec:p2} we examine the case where $p=2$.

\section{Odd primes}\label{sec:odd}
When $p$ is an odd prime and $\rk _p(G)=2$, we show that if the sufficiency condition of 
Proposition \ref{prop:scs} does not hold, then 
$\textrm{Qd}(p)$ must be $p'$-involved in $G$.\\
\\
Before we list the next theorem, we need two definitions as well as a lemma.
\begin{definition}
Let $P \subseteq G$ be a $p$-subgroup. $P$ is said to be \emph{$p$-centric} 
if $Z(P)$ is a Sylow $p$-subgroup of $C_{G}(P)$. Equivalently, $P$ is $p$-centric if and only if 
$C_G(P)=Z(P) \times C' _G(P)$ where 
$C' _G(P)$ is a $p' $-subgroup.
\end{definition}
\begin{definition}[see \cite{g:hlsc}]
Let $P \subseteq G$ be a $p$-centric subgroup. $P$ is said to be \emph{principal $p$-radical}
 if $O_p(N_{G}(P)/P C_{G}(P))=\{1\}$.
\end{definition}

\begin{lemma}[{Diaz, Ruiz, Viruel \cite[Section 3]{drv:plfgr2}}]\label{lem:meta}
Let $G$ be a finite group and $p>2$ a prime with $\rk _p(G)=2$. Let $G_p \in \Syl _p (G)$ 
and $P$ be a principal $p$-radical subgroup of $G$. If $P\neq G_p$ and if $P$ is metacyclic, then 
$P\cong \mathbb{Z}/_{p^n} \times \mathbb{Z}/_{p^n}$ with either $p=3$ or $n=1$.
\end{lemma}

The lemma above follows from Corollaries 3.7 and 3.8 as well as Proposition 3.12 of \cite{drv:plfgr2} .
The argument in \cite{drv:plfgr2} is written in the language of fusion systems, but it can be applied 
to finite groups, which gives Lemma \ref{lem:meta}. We also make use of the following well known lemma, 
which appears in \cite{drv:plfgr2}.

\begin{lemma}[{\cite[Lemma 3.11]{drv:plfgr2}}]\label{lem:pred}
Let $G$ be a subgroup of $\textrm{GL}_2(p)$ for an odd prime $p$. If $O_p(G)=\{1\}$ and if $p$ divides the order of $G$, then 
$\textrm{SL}_2(p)\subseteq G$.
\end{lemma}

\begin{lemma}\label{lem:const}
Let $H$ be a finite group and $p$ an odd prime such that $O_{p'}(H)=\{1\}$ and $C_H(O_p(H))\subseteq O_p(H)$. If $O_p(H)$ is homocyclic abelian 
of rank two and is not a Sylow $p$-subgroup of $H$, then $H$ has a subgroup isomorphic to $\mathrm{Qd}(p)$.
\end{lemma}

\begin{proof}
Let $P=O_p(H) = \mathbb{Z}/_{p^n} \times \mathbb{Z}/_{p^n}$. 
We may consider the group of automorphisms $\aut(\mathbb{Z}/_{p^n} \times \mathbb{Z}/_{p^n})$ as 
two by two matrices with coefficients in $\mathbb{Z}/_{p^n} $ and with  determinants not divisible by $p$. 
Reduction modulo $p$ gives the following short exact sequence of groups:
$$1\rightarrow Q \rightarrow \aut(\mathbb{Z}/_{p^n} \times \mathbb{Z}/_{p^n})\stackrel{\rho }{ \rightarrow } 
\textrm{GL}_2(\mathbb{F}_p) \rightarrow 1$$
with $Q$ a $p$-group. Notice that $H/P \subseteq \aut(\mathbb{Z}/_{p^n} \times \mathbb{Z}/_{p^n})$. 
$Q\cap H/P=\{1\}$ because $O_{p}(H)=P$; therefore, $\rho$ is an injection when 
restricted to $H/P$. Notice that $\rho (H/P)$ has order divisible by $p$, due to $P$ not being a Sylow $p$-subgroup of $H$; 
thus, by Lemma \ref{lem:pred}, $\textrm{SL}_2(p)\subseteq \rho (H/P)$. 
A restriction of $\rho $ is now an isomorphism between a subgroup of $K$ of $H/P$ and $SL_2(\mathbb{F}_p)$.
Let $\alpha \in K$ be preimage  under $\rho $ of the scalar matrix with $-1$ diagonal entries in $SL_2(\mathbb{F}_p)$.
Let $\pi $ be the quotient map $H \rightarrow H/P$. We see that $\alpha $ has a 
preimage under $\pi $ that has order 2. Call this preimage $\beta$. 
$P$ can be viewed as two-dimensional vector space over $\mathbb{F}_{p^n}$ where $p^n$ is the exponent of $P$. 
Notice that $\beta $ acts on elements of $P$ by scalar multiplication by $-1$ so it has no nontrivial fixed points in $P$. 
Let $C$ be the centralizer of $\beta $ in  $H$. 
$\pi |_{C}$ is an isomorphism since $C\cap P=\{1\}$. Also $K\subseteq \pi (C) $ because $K$ is contained in the centralizer of 
$\alpha = \pi (\beta ) $ in $H/P$. It follows that $C$ has a subgroup $L$ isomorphic to $K$, hence, isomorphic to 
$SL_2(\mathbb{F}_p)$. $\langle L , \Omega _1(\bar{P})\rangle $ is a subgroup of $H$, 
which is isomorphic to $\mathrm{Qd}(p)$.
\end{proof}

\begin{remark}
Notice that the hypothesis of Lemma \ref{lem:const} which states that $O_{p'}(H)=\{1\}$ and $C_H(O_p(H))\subseteq O_p(H)$ implies that 
$H$ is $p$-constrained. 
\end{remark}

\begin{theorem}\label{thm:podd}
Let $G$ be a finite group and $p>2$ a prime with $\rk _p(G)=2$. Let $G_p \in \Syl _p (G)$. 
If $\Omega _1 (Z(G_p))$ is not strongly closed in $G_p$ with respect to $G$, then $\mathrm{Qd}(p)$ is $p'$-involved in $G$.
\end{theorem}

\begin{proof} Notice that $\rk (Z(G_p))=1$ by Corollary \ref{cor:omega}. Let $Z=\Omega _1 (Z(G_p))$.  
By hypothesis, $Z$ is not strongly closed in $G_p$ with respect to $G$. Recall that 
Alperin's Fusion Theory \cite{a:siaf} states that given a weak conjugation family $\mathcal{F}$, there 
exists an $(S,T) \in \mathcal{F}$ with $Z\subseteq S \subseteq G_p$ and with $Z$ not strongly closed in 
$S$ with respect to $T$.  The collection of all pairs $(P , N_G(P))$, where $P\subseteq G_p$ is a principal 
$p$-radical subgroup of $G$, is the weak conjugation family of Goldschmidt {\cite[Theorem 3.4]{g:acffg}} 
(see also \cite{j:qsfg}). There exists $P\subset G_p$ that is a principal $p$-radical subgroup of $G$ with $Z\subseteq P$ 
and with $Z$ not strongly closed in $P$ with respect to $N_G(P)$. Notice that $Z\neq P\neq G_p$ and 
that $\rk _p (Z(P))=\rk _p (P)=2$, in particular $p| [N_G(P):P]$. Using the classification of rank two $p$-groups 
for odd primes $p$ by Blackburn 
\cite{b:gcet} (see also \cite{drv:plfgr2,dp:shtr2}), we notice that $P$ must be metacyclic since $\rk _p (Z(P))=2$. 
By Lemma \ref{lem:meta}, $P\cong \mathbb{Z}/_{p^n} \times \mathbb{Z}/_{p^n}$ with either $p=3$ or $n=1$. 
Now consider the group $H=N_G (P)/ C' _G(P)$. We want to show that $H$ satisfies the conditions of Lemma 
\ref{lem:const}. Let $\bar{P}=(P C'_G(P))/C'_G(P)$. 
First $O_p(H)=\bar{P}$ because $P$ is principal $p$-radical. Next let $K=O_{p'}(H)$. $K$ and ${P}$ are both 
normal in $H$ with trivial intersection, so $K\subseteq C_H({P})$ and $K=\{1\}$. $\bar{P}$ 
is not a Sylow $p$-subgroup of $H$ because $P\neq G_p$. Since $P$ is $p$-centric and no $p'$ elements of $H$ commute with $P$,  
$C_H(P)\subseteq P$. Now applying Lemma \ref{lem:const}, we can state that $H$ has a subgroup $Q$ isomorphic 
to $\textrm{Qd}(p)$. $Q$ is thus a subquotient of $G$ obtained by taking a $p'$ quotient, and $\mathrm{Qd}(p)$ is 
$p'$-involved in $G$.
\end{proof}

We point out here some connections between the results of Theorem \ref{thm:podd} and previous work 
by Adem and Smith \cite{as:pcga}. Notice that Theorem \ref{thm:podd} and Proposition \ref{prop:scs} imply 
that any rank 2 group $G$ that does not
have $\mathrm{Qd}(p)$ $p'$-involved in $G$, has a $p$-effective character. 
A $p$-group $P$ is called a Swan group if for any group $G$ containing 
$P$ as a Sylow $p$-subgroup, the mod-p cohomology ring $H^{*}(G)$ is equal to $H^{*}(N_G(P))\cong 
H^{*}(P)^{N_G(P)}$. A group with a Sylow $p$-subgroup that is a Swan group has been shown 
to have a $p$-effective character by Adem and Smith \cite[Section 6]{as:pcga}. Dietz and Glauberman 
(see \cite{mp:chsg}) have shown that any odd order metacyclic group is a Swan group. In addition, 
Diaz, Ruiz and Viruel \cite{drv:plfgr2} have shown that if $P$ is a rank 2 $p$-group for $p>3$, $P$ is Swan unless it is isomorphic to 
a Sylow $p$-subgroup of $\textrm{Qd}(p)$. They also demonstrate that some rank two 3-groups of maximal class 
are Swan.

\begin{remark}\label{rem:glaub}
The results this section are connected to the work of Glauberman on characteristic subgroups and $p$-stability (see \cite{g:psg,g:scps}). To discribe the connection, we first need to define the Thompson subgroup of a $p$-group $P$: the Thompson subgroup $J(P)$ is then the subgroup of $P$ generated by the abelian subgroups of $P$ that have maximal order among all abelian subgroups of $P$. Glauberman in \cite{g:psg} shows that if $G$ is a finite group and $p$ an odd prime with $G_p\in \Syl _p(G)$ such that $\mathrm{Qd}(p)$ is not involved in $G$, then $Z(J(G_p))$ is strongly closed in $G_p$ with respect to $G$. If, in addition, $Z(J(G_p))$ is contained in every elementary abelian $p$-subgroup of $G_p$ with the same rank as $G$, then $G$ has a $p$-effective character; in particular, $G$ has a $p$-effective character for an odd prime $p$ if $G$ is $\mathrm{Qd}(p)$-free and $J(G_p)=G_p$. Notice that this result does not depend on the rank of $G$. 
\end{remark}

\section{Prime two}\label{sec:p2}
The argument in this section will be done in two parts. 
First we will show, regarding the prime two, that if the sufficiency condition of Proposition \ref{prop:scs} 
does not hold for $G$ and for prime $2$, 
then $G_2\in \Syl _2(G)$ is either dihedral, semi-dihedral, or wreathed. In the second part we will 
deal with these exceptional cases. We start with two definitions.

\begin{definition}[{\cite[p. 191]{g:fg}}]
Recall that a 2-group is semi-dihedral (sometimes called quasi-dihedral) if it is 
generated by two elements $x$ and $y$ subject to the relations that $y^{2}~=~x^{2^{n}}=1$ and 
$yxy^{-1}=x^{-1+2^{n-1}}$ for some $n\geq 3$.
\end{definition}
\begin{definition}[{\cite[p. 486]{g:fg}}]\label{def:wr}
A 2-group is called wreathed if it is generated by three elements 
$x$, $y$, and $z$ subject to the relations that $x^{2^{n}}=y^{2^{n}}=z^{2}=1$, $xy=yx$, and 
$zxz^{-1}=y$ with $n \geq 2$.
\end{definition}

Notice that the dihedral group of order 8 is wreathed. 
Now that we have these definitions, we will state Proposition \ref{prop:p2}, which 
is crucial in our discussion of 2-subgroups of rank two. The proof of Proposition 
\ref{prop:p2} follows the proof of a theorem by Alperin, Brauer, and Gorenstein 
\cite[Proposition 7.1]{abg:fsg2r2}.

\begin{proposition}[See Proposition 7.1 of \cite{abg:fsg2r2} and its proof.]\label{prop:p2}
Let $G$ be a finite group and $G_2\in \Syl _2(G)$. Suppose that $\rk (G_2)=2$. If 
$\Omega _1 (Z(G_2))$ is not strongly closed in $G_2$ with respect to $G$, then $G_2$ is either 
dihedral, semi-dihedral, or wreathed.
\end{proposition}

Below we briefly review the proof of Proposition 7.1 of \cite{abg:fsg2r2}, including the changes 
needed to prove Proposition \ref{prop:p2}.

Again notice that $\rk (Z(G_2))=1$ by Corollary \ref{cor:omega}. 
If $G_2$ does not contain a normal rank 2 elementary abelian subgroup, then $G_2$ 
must be either dihedral or semi-dihedral (see \cite[Theorem 5.4.10]{g:fg}). 
We may assume that $G_2$ has a normal rank 2 elementary abelian subgroup, which we will call $V$. 
Letting $T=C_{G_2}(V)$, we notice that $[G_2:T]\leq 2$. Also $G_2 \neq T$ since $\rk (Z(G_2))=1$; 
therefore, $[G_2:T]= 2$.\\
\\
Letting $A_G(V)=N_G(V)/C_G(V)$, we notice that $A_G(V)$ is isomorphic to a subgroup of 
$\aut (V)\cong \Sigma _3$ and that $2\bigr| |A_G(V)| \bigl.$. Thus either $A_G(V)\cong \Sigma _3$ or $2=|A_G(V)|$.\\
\\
From the argument in \cite{abg:fsg2r2}, we get the following results: if $A_G(V)\cong \Sigma _n$, then 
$G_2$ is wreathed, including the case of $D_8$; and if $|A_G(V)|=2$, then $G_2 \cong D_8$.
This finishes the discussion of Proposition \ref{prop:p2}.\\
\\
We have seen that if $G$ is a rank two finite group, it has a 
2-effective character unless a Sylow 2-subgroup is dihedral, semi-dihedral, or wreathed. 
In the rest of this section, we will show that in each of these three cases $G$ also has a 2-effective 
character, which will complete the proof that a rank two finite group has a 2-effective 
character. The author proved the following three lemmas in his doctoral thesis 
\cite{j:vbobg}; the proofs are included here for completeness. We will start with a 
lemma that will be used in the proof of Proposition \ref{prop:p2ex}.

\begin{lemma}\label{lem:ann}
Let $G$ be a finite group, let $n=\rk(G)$, let $p$ a prime divisor of $|G|$ with $\rk _{p}(G)=n$, 
and let $G_p\in \Syl _p (G)$. If $G_p$ is abelian, $G_p \lhd G$, 
or $G$ is $p$-nilpotent, then $G$ has a $p$-effective character.
\end{lemma}
\begin{proof} In each case we will show that $G$ has a $p$-effective character by 
showing that $Z(G_p)$ is strongly closed in $G_p$ with respect to $G$ using Proposition 
\ref{prop:scs}. If $G_p$ is abelian, $Z(G_p)=G_p$ is obviously strongly closed in $G_p$ with 
respect to $G$ (see Corollary \ref{cor:omega}). In the case where $G_p\lhd G$, $Z(G_p) \lhd G$ 
because $Z(G_p)$ is a characteristic 
subgroup of $G_p$; therefore, $Z(G_p)$ is strongly closed in $G_p$ with respect to $G$. If $G$ 
is $p$-nilpotent, then two elements of $G_p$ are conjugate in $G$ if and only if they are conjugate 
in $G_p$. Thus, in this case as well, $Z(G_p)$ is strongly closed in $G_p$ with respect to $G$. 
\end{proof}

\begin{lemma}\label{lem:dih}
If $P$ is a dihedral or semi-dihedral 2-group such that $|P|=2^n$ with $n\geq 3$, then there is a 
character $\chi $ of $P$ such that 
\begin{displaymath}
\chi (g)=\left\{\begin{array}{ll}
3*2^{n-3} & \textrm{if $g=1$}\\
-2^{n-3} & \textrm{if $g$ is an involution}\\
2^{n-3} & \textrm{otherwise.}
\end{array} \right. 
\end{displaymath}
\end{lemma}

\begin{proof} Let $N$ be the commutator subgroup 
$[P,P]$ and notice that $P/N \cong \mathbb{Z}/2 \times \mathbb{Z}/2$. Since $Z(P)\cong \mathbb{Z}/2$, $Z(P)\subseteq N$; therefore, there is a non-trivial irreducible 
character $\lambda  $ of $P/N$ with $\lambda  (gN)=1 $ for $g\in Z(P)$   such that if $g\in P \setminus Z(P)$ is an involution, then 
$\lambda (gN)=-1$. Let $\varphi (g)= \lambda (gN)$, which is a character of $P$.\\
\\
We will define another character of $P$ by first inductively defining a character of each $D_{2^m}$. 
We start by letting $\psi _3$ be the character of $D_8$ such that 
\begin{displaymath}
\psi _{3}(g)=\left\{\begin{array}{ll}
2 & \textrm{if $g=1$}\\
-2 & \textrm{if $g \in Z(D_8)\setminus \{1\}$}\\
0 & \textrm{otherwise.}
\end{array} \right.
\end{displaymath}
For the induction, let $\psi _{k}= \ind _{D_{2^{k-1}}}^{D_{2^k}} \psi _{k-1}$. 
$\psi _k$ is a character of $D_{2^k}$ for each $k\geq 3$. We see that 
if $P$ is dihedral, then $\chi = 2^{n-3} \varphi + \psi _n$; therefore, Lemma \ref{lem:dih} holds for dihedral 2-groups.
On the other hand, if $P$ is semidihedral, $\chi = 2^{n-3} \varphi + \ind _{D_{2^{n-1}}}^{P} \psi _{n-1}$.
\end{proof}

\begin{lemma}\label{lem:wr}
Let $G$ be a finite group with $S\in \Syl _2 (G)$ a wreathed 2-group. If $\rk(G)=2$ and $G$ has no normal subgroup of index 2, then $G$ has a 2-effective character.
\end{lemma}
\begin{proof} Let $S$ be a Sylow 2-subgroup of $G$. We will say that $S$ is 
generated by $x$, $y$, and $z$ as in Definition \ref{def:wr}. 
Let $\alpha $ be any primitive $(2^{n})^{th}$ root of 
unity. We will define a degree 3 complex representation of $S$ as follows.
$ \kappa :S \rightarrow GL_{3}(\mathbb{C}) $
is given by
\[ \kappa (x) = \left( \begin{array}{c c c} \alpha  & 0 & 0 \\ 0 & \alpha  & 0 \\
0 & 0 & \alpha ^{-2} \end{array} \right) \]
\[ \kappa (y) = \left( \begin{array}{c c c} \alpha & 0 & 0 \\ 0 & \alpha ^{-2} & 0 \\
0 & 0 & \alpha \end{array} \right) \]
\[ \kappa (z) = \left( \begin{array}{c c c } -1 & 0 & 0 \\ 0& 0 & 1 \\ 0& 1 & 0 \end{array} \right) .\]
Let $\nu $ be the complex character of $S$ associated to $\kappa$: for each $g\in S$, $\nu (g)= 
\tr (\kappa (g))$. We will now show that $\nu $ is a 2-effective character 
for $G$.\\
\\
First notice the character values on involutions. The elements of $S$ that are involutions are 
$x^{2^{n-1}}$, $y^{2^{n-1}}$, $(xy)^{2^{n-1}}$, and those conjugate in $S$ to $z$. Clearly $\nu $ 
takes the value $-1$ for each of these elements. If $H$ is a rank 2 elementary abelian subgroup of $S$, then 
$[\nu |_{H},1_{H}]=0$.\\
\\
Let $S_{0}$ be the subgroup of $S$ generated by $x$ and $y$. $S_{0}$ is 
an abelian normal subgroup of index 2. We emphasize that $z$ is conjugate to an element of $S_0$, 
in particular to $x^{2^{n-1}}$, and $[N_G(S_0):C_G(S_0)]=6$ (see \cite[Proposition 2G*]{bw:spfgwss}). Now we look at fusion of elements of $S_{0}$. 
Any two elements of $S_{0}$ that are conjugate in $G$ are conjugate in $N_{G}(S_{0})$ (see \cite[Proposition 2D]{bw:spfgwss}). 
$N_{G}(S_{0})$ is generated by $C_{G}(S_{0})$, the element $z$, and an element  $\eta \in G$ such that $\eta x\eta ^{-1}=
y$,  $\eta y\eta ^{-1}=x^{-1}y^{-1}$, and $\eta (x^{-1}y^{-1})\eta ^{-1}=x$ \cite[Proposition 2B$_{\textrm{I}}$]{bw:spfgwss}.\\
\\
To show that $\nu $ respects fusion of elements in $S_0$, it is enough to show 
that if $g \in S_{0}$, then $\nu (\mu )=\nu (z g z^{-1})=\nu (\eta g \eta ^{-1})$, assuming that such an $\eta $ exists.
Notice that for any $g \in S_{0}$, there are integers $k$ and $l$ with $g =x^{k}y^{l}$. 
We see then that $z g z^{-1}=x^{l}y^{k}$ and $ \nu (x^{k}y^{l})=\alpha ^{k+l}+\alpha ^{k-2l} +\alpha ^{-2k+l}=\nu _{i} (x^{l}y^{k})$. Also $\eta g \eta ^{-1}=x^{-k}y^{l-k}$ and 
$ \nu (x^{k}y^{l})=\alpha ^{k+l}+\alpha ^{k-2l}+\alpha ^{-2k+l}= \alpha ^{l-k-k}+
\alpha ^{-k-2(l-k)}+\alpha ^{2k+l-k}=\nu ( x^{-k}y^{l-k})$.\\
\\
Next we will check fusion in elements of order at most $2^n$. If $g \in S \setminus 
S_{0}$ is of order at most $2^{n}$, then there is an integer $h$ with $g $ conjugate to 
$x^{-h}y^{-h}z$ in $S$ (see \cite[Lemma 4]{bw:spfgwss}). In particular, 
$\nu (g )= \nu  (x^{-h}y^{-h}z)$. Recall that $x^{-h}y^{-h}z$ is conjugate to 
$x^{2^{n-1}-h}y^{-h}$ in $G$ (see \cite[Proposition 2F$_{\textrm{a}}$]{bw:spfgwss}). We see that
$\nu  (x^{-h}y^{-h}z)=-\alpha ^{-2h} = \alpha ^{2^{n-1}} \alpha ^{-2h} =
\nu  (x^{2^{n-1}-h}y^{-h})$.\\
\\
At this point we have finished the fusion of elements of order at most $2^m$. Notice that all elements of order larger than $2^m$ must have order $2^{m+1}$ and must be in $S\setminus S_0$. Suppose that $g\in S$ is of order $2^{m+1}$ and is conjugate to another element $g^{\prime }$ of $S$. It immediately follows that $g, g^{\prime} \in S \setminus S_0$ and that $g$ and $g^{\prime }$ are conjugate in $S$ (see \cite[Proposition 2E]{bw:spfgwss}). So any character of $S$, in particular $\nu $, must agree on $g$ and $g^{\prime}$. 
This concludes the proof of Lemma \ref{lem:wr}.
\end{proof}

\begin{proposition}[Jackson \cite{j:vbobg}] \label{prop:p2ex}
If $G$ is a finite group with a dihedral, semi-dihedral, or wreathed Sylow 2-subgroup such that 
$\rk(G)=2$, then $G$ has a 2-effective character.
\end{proposition}

\begin{proof} Let $G_2\in \Syl _2(G)$ and let $n$ be the integer with $|G_2|=2^n$. If $G_2$ is either 
dihedral or semi-dihedral, then there is a character of $G_2$ given in Lemma \ref{lem:dih} such that 

\begin{displaymath}
\chi (g)=\left\{\begin{array}{ll}
3*2^{n-3} & \textrm{if $g=1$}\\
-2^{n-3} & \textrm{if $g$ is an involution}\\
2^{n-3} & \textrm{otherwise     .}
\end{array} \right. 
\end{displaymath}

Since $\chi $ is constant on elements of the same order, it must respect fusion in $G$ as in Remark 
\ref{rem:peff}. 
It is also clear that if $E\subseteq G_2$ is a rank 2 elementary abelian subgroup, 
$[\chi |_{E},1_{E}]=0$ (see Remark \ref{rmk:reg}). $\chi $ is, thus, a 2-effective character of $G$.\\
\\
We may now assume that $G_2$ is wreathed and is generated by $x$, $y$, and $z$ as 
in Definition \ref{def:wr}. In addition, let $S_0$ be the abelian subgroup of $G_2$ generated by $x$ and $y$. Since $G_2$ is wreathed, one of the following holds 
(see \cite[Proposition 2, p. 11]{abg:fgss}):
\begin{enumerate}
\item $G$ has no normal subgroups of index 2,
\item $Z(G_2)$ is strongly closed in $G_2$ with respect to $G$,
\item $G$ has a normal subgroup $K$ with Sylow 2-subgroup $S_{0}$ generated by $x$ and 
$y$, or
\item $G$ is 2-nilpotent.
\end{enumerate}

Case 1 was treated in Lemma \ref{lem:wr}. Case 2 has been treated in Proposition 
\ref{prop:scs}.\\
\\
We now treat case 3. Suppose $G$ has a normal subgroup $K$ with Sylow 2-subgroup $S_{0}$ generated by $x$ and 
$y$. Notice that $G_2\cap K=S_0$. We see that $S_0$ is weakly closed 
in $G_2$ with respect to $G$. By a standard result of Burnside (see \cite[37.6]{a:fgt}), 
$N_G(S_0)$ controls fusion in $C_{G_2}(S_0)=S_0$.  
Also notice that any element of $G_2\setminus S_0$ is conjugate in $G$ to only those elements to which
it is conjugate in $G_2$. Using the Frattini factor subgroup of $S_0$, we notice that $N_G(S_0)/C_G(S_0)$ is 
isomorphic to a subgroup of $\Sigma _3$ (see \cite[Page 264]{bw:spfgwss} or \cite[Page 12]{abg:fgss}). 
We may assume that $Z(G_2)$ is not strongly closed in $G_2$ with respect to $G$, otherwise the group is treated 
in case 2. In particular, this implies that $\Omega _1(Z(G_2))=\langle (xy)^{2^{n-1}}\rangle $ is not strongly closed in 
$G_2$ with respect to $G$. The involution $(xy)^{2^{n-1}}$ must then be conjugate in $N_G(S_0)$ to another involution in 
$S_0$. Thus, $N_G(S_0)/C_G(S_0)$ must not be a two group and must also have order divisible by two since $S_0 \lhd G_2$. 
$N_G(S_0)/C_G(S_0)\cong \Sigma _3$, and we notice this symmetric group 
permutes the elements 
$x$, $y$, and $xy$; furthermore, 
the character $\chi $ of $G_2$ given in Lemma \ref{lem:wr} is again a 2-effective character in this case.\\
\\
Case 4 was treated in Lemma \ref{lem:ann}.
\end{proof}

To conclude Section \ref{sec:p2}, we state the following theorem whose proof has now been completed.

\begin{theorem}\label{thm:2eff}
Every rank two finite group $G$ has a 2-effective character.
\end{theorem}

\section{Counterexample}\label{sect:count}
Looking at Section \ref{sec:p2}, one may wonder if a similar argument can be made for odd primes.  
However, this is not the case as we will show in the following proposition.

\begin{proposition}
Let $p$ be an odd prime and $G$ be a rank two finite group. If $\textrm{Qd}(p)$ is involved in $G$ and if $G$ and $\textrm{Qd}(p)$ 
have isomorphic sylow $p$-subgroups, then $G$ does not have a $p$-effective character. In particular, $\textrm{Qd}(p)$ does not have a $p$-effective character for an odd prime $p$.
\end{proposition}

\begin{proof} Let $G_p \in \Syl _p (G)$. Notice that $G_p$, which is isomorphic to a Sylow $p$-subgroup of $\textrm{Qd}(p)$, is an 
extra-special $p$-group of size $p^3$ and exponent $p$. 
Let $\chi _{1}, \dots , \chi _{n}$ be the irreducible characters of $G_p$ and assume that $\chi _{1}$ is the 
trivial character. Notice from the structure of $G_p$ that if $\chi _{i}(1)\neq 1$, then for $g\in G_p$, 
$\chi _{i}(g)\neq 0$ if and only if $g\in Z(G_p)$. Also notice that $Z(G_p)$ is cyclic of order $p$ and 
$Z(G_p)$ is not strongly closed in $G_p$ with respect to $G$. Let $x\in Z(G_p)$ and $y\in G_p \setminus 
Z(G_p)$ such that there is a $g\in G$ with $y^g=x$. Suppose 
$\chi $ is a character of $G_p$ that is a $p$-effective character of $G$. There exists $a_{1}, \dots 
, a_{n} \in \mathbb{Z}_{\geq 0}$ such that $\chi = \sum_{i=1}^{n} a_{i}\chi _{i}$. Since $\chi $ is 
a $p$-effective character of $G$, it respects fusion in $G$, which implies $\chi (x) = \chi (y)$ and 
$a_{1} =0$. Suppose that $a_{i}>0$ for some $i>1$ with $\chi _{i}(1)=1$. Fixing this $i$, let 
$E_{i}=\{g\in G_p | \chi _{i}(g)=1\}$. Notice 
that for such an $i$, $E_{i}$ is a rank two elementary abelian subgroup of $G_p$. This implies that 
$[\chi _{i} |_{E_{i}},1_{E_{i}}]=1$, so $[\chi |_{E_{i}},1_{E_{i}}]>0$, which contradicts the assumption that $\chi $ 
is a $p$-effective character of $G$; therefore, 
for each $i$ such that $\chi _{i}(1)=1$, $a_{i}=0$. 
Now \[ \chi (x) = 
\sum_{i \textrm{\footnotesize { such that }}\chi _{i}(1)\neq 1} a_{i} \chi _{i} (x) + 
\sum_{i \textrm{\footnotesize { such that }}\chi _{i}(1)= 1} a_{i}\]
\[= 
\sum_{i \textrm{\footnotesize { such that }}\chi _{i}(1)= 1} a_{i} \chi _{i}(y) = \chi (y).\] 
$\chi (y)=0$; therefore, $\chi (x)=0$.  
$\chi (z)=0$ for all $z\in Z(G_p)\setminus \{1\}$ because $Z(G_p)$ is a cyclic group of order $p$. 
Notice that in showing $\chi (z)=0$ for each 
$z\in Z(G_p)\setminus \{1\}$, we see that $\chi (g)=0$ for all $g\in G_p\setminus Z(G_p)$ 
by the structure of $G$; thus, $\chi (g)=0$ for all 
$g\in G_p\setminus \{1\}$. $\chi $ must be identically zero, which contradicts 
the definition of $p$-effective character.
\end{proof}

\begin{lemma}\label{lem:sub}
Let $G$ be a finite group with $p$ a prime dividing $|G|$ and $H\subseteq G$. Suppose 
that $p$ divides $|H|$ and that $\rk _p(G)=\rk _p(H)$. If $G$ has 
a $p$-effective character, so does $H$.
\end{lemma} 

\begin{proof} We may assume that $\rk _p (H)=\rk (G)$. 
Let $G_p\in \Syl _p (G)$ such that $H_p=G_p\cap H \in \Syl _p (H)$. Let $\chi $ be a character of 
$G_p$ that is a $p$-effective character of $G$. $\chi |_{H_p}$ is a character of 
$H_p$, which is not identically zero. $\chi |_{H_p}$ respects fusion in $H$
because $\chi $ respects fusion in $G$. Any maximal rank elementary abelian subgroup of $H_p$ is 
also a maximal rank elementary abelian subgroup of $G_p$; therefore, $\chi |_{H_p}$ is a $p$-effective character of $H$.
\end{proof}

Combining Lemma \ref{lem:sub} with Theorems \ref{thm:podd} and \ref{thm:2eff}, 
we get the following theorem from which Remark \ref{rem:psl3} follows:

\begin{theorem}\label{thm:iff}
Let $G$ be a finite group of rank two and let $p$ be a prime dividing $|G|$. $G$ has $p$-effective 
character if and only if either $p=2$ or both $p>2$ and $G$ does not $p'$-involve $\textrm{Qd}(p)$.
\end{theorem}

In showing that a finite group $G$ of rank two acts freely on a finite complex 
$Y\simeq \mathbb{S}^{n} \times \mathbb{S}^m$, we show that $G$ acts on a finite 
complex $X\simeq \mathbb{S}^m$ with isotropy groups of rank one and then we apply Theorem 
\ref{thm:as1}. Ozgun Unlu \cite{u:dd} and Grodal \cite{g:pc} have shown separately that for each odd prime 
$p$, $\textrm{Qd}(p)$ cannot act on any finite complex homotopy equivalent to a sphere with rank one 
isotropy groups. From both of their proofs it is clear that for any odd prime p and rank two group $G$ that $p'$-involves $\textrm{Qd}(p)$, $G$ cannot act on 
any finite complex homotopy equivalent to a sphere with rank one isotropy groups.\\
\\
Combining this discussion with Theorem \ref{thm:iff} and the discussion in 
Section \ref{sec:pre}, we conclude with the following proposition:

\begin{proposition}
Let $G$ be a finite group of rank two. $G$ acts on some finite complex homotopy equivalent to a 
sphere with rank one isotropy groups if and only if for each odd prime $p$, $G$ does not $p'$-involve $\textrm{Qd}(p)$.
\end{proposition}

\end{document}